\documentclass{amsart}
\usepackage{amssymb}
\usepackage{a4wide}
\usepackage[english]{babel}
\newcommand{\Ps}{\mathbf{P}}
\newcommand{\Z}{\mathbf{Z}}

\newcommand{\Q}{\mathbf{Q}}\newcommand{\M}{\mathcal{M}}

\newcommand{\str}{\mathcal{O}}
\newcommand{\ra}{\rightarrow}

\renewcommand{\phi}{\varphi}
    \newtheorem{Lem}{Lemma}[section]
    \newtheorem{Prop}[Lem]{Proposition}
    \newtheorem{Thm}[Lem]{Theorem}
    
    \newtheorem{Cor}[Lem]{Corollary}
    
    \newtheorem{Ass}[Lem]{Assumption}
    \newtheorem{Not}[Lem]{Notation}
   \theoremstyle{definition}
    \newtheorem{Def}[Lem]{Definition}
    
    \newtheorem{Rem}[Lem]{Remark}
    \DeclareMathOperator{\rank}{rank}

\begin{document}
\title[Classification of fibrations]{Classification of all Jacobian elliptic fibrations on certain $K3$ surfaces}\label{classi}
\author{Remke Kloosterman}
\date{\today}
\thanks{The author would like to thank Marius van der Put, Jaap Top and Joseph Steenbrink for the many comments on an earlier version of this article. This article is based on a chapter in the author's PhD thesis \cite[Chapter 5]{proefschrift}.}
\begin{abstract}
In this paper we classify all configurations of singular fibers of elliptic fibrations on the double cover of $\Ps^2$ ramified along six lines in general position.
\end{abstract}
\maketitle

\section{Introduction}
In this article we classify all possible bad fiber configurations on
Jacobian elliptic fibrations on the $K3$ surface $X$, which is the minimal model  of the double
cover of $\Ps^2$ ramified along six lines in `general position'. When we say that six lines are
in `general position' we mean that the rank of the N\'eron-Severi 
group $NS(X)$ is $16$.

The strategy we use is purely geometric, and very similar to Oguiso's
classification of Jacobian elliptic fibrations on the Kummer surface
of the product of two non-isogenous elliptic curves (\cite{Ogui}). It
seems possible to extend this method to the case of $K3$ surfaces,
 which are birational to a double cover of $\Ps^2$ ramified along a
 sextic curve, such that the N\'eron-Severi group of the $K3$ surface
 is generated by the (reduced) components of the pull-back of the
 branch divisor together with all divisors obtained by resolving
 singularities on the double cover and the pull-back of a general line
 in $\Ps^2$.

Here we give the full list of possibilities.

\begin{Thm}\label{thminfinite} Let $X$ be as before. Suppose  $\pi : X \ra \Ps^1$ be an
  elliptic fibration with positive Mordell-Weil rank. Then the
  configuration of singular fibers is contained in the following list.
\[\begin{array}{clcc}
\mbox{Class}& \mbox{Configuration of singular fibers}& & MW\mbox{-rank}\\
1.1 & I_{10}\; I_{2}\; aII\;bI_1& 2a+b=12 & 4 \\
1.2 & I_{8} \;I_{4} \;aII\;bI_1& 2a+b=12  &4 \\
1.3 & 2I_{6} \; aII\;bI_1& 2a+b=12 & 4 \\
1.4 & IV^*\; I_4 \; aII\;bI_1& 2a+b=12  &5\\ 
\end{array}\]
Conversely, for each class there exists $a,b$ such that these fibrations occur.
\end{Thm} 
We did not
establish the possible values of $a$ and $b$.  By Proposition~\ref{NLbnd} we know that generically $a\leq 4$ in $1.1, 1.2$ and $1.3$, $a\leq 5$ in $1.4$. We believe that for each of the classes $1.1-1.4$ we have that generically $a=0$.

\begin{Thm}\label{thmfinite} Let $X$ be as before. Let $\pi: X \ra \Ps^1$ be a Jacobian
  elliptic fibration with finite Mordell-Weil group. Then the
  configuration  of singular fibers is contained in Table~\ref{tabpos}.
(In this table we list all fiber types, different from $III$, $I_2$, $II$,
  $I_1$, the structure of $MW(\pi)$ and two quantities, namely $iii+i_2$ and
  $3iii+2i_2+2ii+i_1$, where $iii$ means the number of fibers of type
  $III$ etc.)
Conversely, for each class there exists a Jacobian elliptic fibration $\pi : X \ra \Ps^1$ with the given configuration of singular fibers.
\end{Thm}

\begin{table}
\[\begin{array}{c|l|c||cc||cc}
\mbox{Class}&\mbox{Configuration} & MW(\pi) &  iii+i_2 & 3iii+2i_2+2 ii+i_1\\
2.1 & II^* & 1 & 6 &14\\
2.2 & III^* & \Z/2\Z & 7&15 \\
2.3 & III^* I_0^* &1&   3 & 9 \\
2.4 & I_6^* &  1 & 4 &12 \\
2.5 & I_4^*  & \Z/2\Z & 6&14\\
2.6 & I_4^* I_0^* & 1 & 2&8\\
2.7& I_2^*&(\Z/2\Z)^2 &  8& 16\\
2.8 & I_2^* I_0^*&\Z/2\Z&  4 & 8 \\
2.9 & 2 I_2^* & 1  & 2&8\\
2.10 & I_2^*\; 2I_0^* & 1 & 0 & 8\\
2.11& 2I_0^* & (\Z/2\Z)^2 & 6 &12\\
2.12 & 3I_0^* & \Z/2\Z & 2 & 6\\
\end{array}\]
\caption{List of possible configurations}\label{tabpos}
\end{table}

\begin{Cor}\label{corfin} Let $X$ be as before, and let $\pi: X \ra \Ps^1$ be an
  elliptic fibration with finite Mordell-Weil group.        If the six lines are sufficiently general
  then the configuration of singular fibers of $\pi$ is contained in
  Table~\ref{tabposgen}.
Conversely, for each class there exists a Jacobian elliptic fibration $\pi : X \ra \Ps^1$ with the given configuration of singular fibers.
\end{Cor}
\begin{table}
\[\begin{array}{c|l|c||cc||cc}
\mbox{Class}&\mbox{Configuration} & MW(\pi) & i_2 & i_1 \\
2.1 & II^* & 1 & 6 &2\\
2.2 & III^* & \Z/2\Z & 7& 1 \\
2.3 & III^* I_0^* & 1 & 3 & 3 \\
2.4 & I_6^* &  1 & 4&4 \\
2.5 & I_4^*  & \Z/2\Z & 6 & 2\\
2.6 & I_4^* I_0^* & 1 & 2& 4\\
2.7& I_2^*&(\Z/2\Z)^2 & 8 & 0 \\
2.8 & I_2^* I_0^*&\Z/2\Z& 4 & 0  \\
2.9 & 2 I_2^* & 1 & 2 & 4\\
2.10 & I_2^*\; 2I_0^* & 1 &0 & 8 \\
2.11& 2I_0^* & (\Z/2\Z)^2 & 6 & 0\\
2.12 & 3I_0^* & \Z/2\Z & 2 & 2\\
\end{array}\]
\caption{List of possible configurations on a general $X$}\label{tabposgen}
\end{table}
This corollary is an immediate consequence from
Theorem~\ref{thmfinite} and Proposition~\ref{NLbnd}.

\begin{Rem}Note that in the cases 2.7, 2.8, 2.10 and 2.11, the fiber
configuration in Theorem~\ref{thmfinite} and Corollary~\ref{corfin}
are the same.\end{Rem}

The generality condition is very essential in our strategy: one of the
key tools in this article uses the fact that the involution $\sigma$
on $X$ induced by the double-cover involution acts
trivially on the N\'eron-Severi group $NS(X)$. Our definition of
general position implies that $\sigma$ acts trivial on $NS(X)$.

The organization of this article is as follows. In Section~\ref{Prel} we introduce some basic definitions and notation. In Section~\ref{scur}
we start with recalling some facts about curves on $K3$ surfaces. In
Section~\ref{skod} we recall
some standard facts concerning singular fibers of elliptic fibrations and
some special results on elliptic surfaces. In Section~\ref{sNS}  we study the
N\'eron-Severi group of a double cover of $\Ps^2$, ramified along six
lines in general position. In Section~\ref{sNL} we give a variant of
\cite[Proposition 2.4.6]{NL}, which implies that Corollary~\ref{corfin}
follows from Theorem~\ref{thmfinite}. 
In Section~\ref{sps} we list all types of possible singular fibers for an
elliptic fibration on $X$ and we  count the number of pre-images of the six
lines contained in each fiber type. In Section~\ref{sca} we classify all fibrations in which all special
rational curves  (the pre-images of the six branch lines) are
contained in the singular fibers, thereby proving Theorem~\ref{thminfinite}. In
Section~\ref{scb} we prove Theorem~\ref{thmfinite}.

Every proof of the actual existence of a fibration presented in this
article runs  as follows: First we give
an effective divisor $D$, with $D^2=0$ and such that there exists an
irreducible curve $C\subset X$ with $D.C=1$. It is easy to see that
then $|D|$ defines an elliptic fibration $\pi: X \ra \Ps^1$ (see
Lemma~\ref{fiblem}), that $C\cong \Ps^1$, and that  $(\pi|_C)^{-1}: \Ps^1 \ra X$ is a section.

In this article all fibrations, sections and components of singular fibers are defined over the field of definition of the six lines.
\section{Definitions \and Notation}
\label{Prel}
%

\begin{Def}\label{defbas} An \emph{elliptic surface} is a triple
  $(\pi,X,C)$ with $X$ a surface, $C$ a curve, $\pi$ is a morphism
  $X\rightarrow C$, such that almost all fibers are irreducible genus
  1 curves and $X$ is relatively minimal, i.e., no fiber of $\pi$
  contains an irreducible rational curve $D$ with $D^2=-1$.

  We denote by $j(\pi): C \rightarrow \Ps^1$ the rational function such that $j(\pi)(P)$ equals the $j$-invariant of $\pi^{-1}(P)$, whenever $\pi^{-1}(P)$ is non-singular.

A \emph{Jacobian elliptic surface} is an elliptic surface together with a section $\sigma_0: C \rightarrow X$ to $\pi$.
  The set of sections of $\pi$ is an abelian group, with $\sigma_0$ as the identity element. Denote this group by $MW(\pi)$.

By an {\em elliptic  fibration} on $X$ we mean that we give a surface $X$ a structure of an elliptic surface.

Let $NS(X)$ denote the  group of divisors modulo algebraic equivalence. We call $NS(X)$ the N\'eron-Severi group of $X$. Denote by $\rho(X)$ the rank of $NS(X)$. We call $\rho(X)$ the {\em Picard number} of $X$.
\end{Def}

Recall the following theorem.
\begin{Thm}[{Shioda-Tate (\cite[Theorem 1.3 \& Corollary
    5.3]{Sd})}]\label{ST}  Let $\pi:X\rightarrow C$ be a Jacobian
  elliptic surface, such that at least one of the fibers of $\pi$ is singular. Then the N\'eron-Severi
  group of $X$ is generated by the classes of $\sigma_0(C)$, a
  non-singular fiber, the components of the singular fibers not
  intersecting $\sigma_0(C)$, and the generators of the Mordell-Weil
  group. Moreover, let $S$ be the set of points $P$ such that $\pi^{-1}(P)$ is singular. Let $m(P)$ be the number of irreducible components of $\pi^{-1}(P)$, then
  \[ \rho(X) =2+ \sum_{P \in S} (m(P)-1)+\rank(MW(\pi)) .\]
\end{Thm}


\section{Curves on $K3$ surfaces}\label{scur}
In this section we give some elementary results on curves on $K3$ surfaces.
\begin{Lem}\label{selfint} Suppose $D$ be a smooth curve on a $K3$ surface $X$. 
Then
\[ g(D)=1+\frac{D^2}{2}\]\end{Lem}
\begin{proof}
Since the canonical bundle $K_X$ is trivial, the adjunction formula
for a divisor on a $K3$ surface is  $ 2p_a(D)-2= D^2$ (see
\cite[Proposition V.1.5]{Har}). Since $p_a(D)=g(D)$ for a smooth
curve, this implies the result.
\end{proof}

\begin{Lem}\label{fiblem} Suppose $D$ is an effective divisor on a $K3$ surface $X$ with $p_a(D)=1$. Then
  $|D|$ defines an elliptic fibration $\pi: X \ra \Ps^1$.
Every effective connected divisor $D'$ such that $D\cdot D'=D'^2=0$ is a fiber of $\pi$.
\end{Lem}

\begin{proof} From the  the adjunction formula
\cite[Proposition V.1.5]{Har} it follows that $D^2=0$. 

 Applying Riemann-Roch \cite[Theorem V.1.6]{Har} yields
\[ \dim H^0(X,\str_X(D)) +\dim H^0(X,\str_X(-D)) \geq 2. \]
Since $\dim H^0(X,\str_X(-D))=0$ ($D$ is effective), we obtain that
$\dim H^0(X,\str_X(D))>1$. This, combined with the fact that $D^2=0$, implies
that $|D|$ is base-point-free. 
So we can apply Bertini's theorem \cite[Theorem II.8.18 \& Remark III.7.91]{Har}, hence there is an
irreducible curve $F$ linearly equivalent to $D$. By
Lemma~\ref{selfint} this curve has genus 1.

The exact sequence
\[ 0\ra \str_X \ra \str_X(F)\ra \str_F \ra 0,\]
together with  the facts  $H^1(X,\str_X)=H^2(X,\str_X(F))=0$ and
$\dim H^2(X,\str_X)=1$, gives that $H^1(X,\str_X(F))=0$, whence by
Riemann-Roch \cite[Theorem V.1.6]{Har} we obtain that the dimension of $
H^0(X,\str_X(F))$ equals 2, so ``the'' morphism associated to $|D|$,
$\pi: X \ra \Ps^1$ is an elliptic fibration.

Since $D'$ is effective and $D'\cdot D=0$, we obtain that no irreducible
component of $D'$ is  intersecting $D$. Since $D'$ is connected we
obtain that $\pi(D')$ is a point. From $D'^2=0$ it follows that
$p_a(D')=p_a(D)$, hence $D'$ is a fiber.
\end{proof}

\section{Kodaira's classification of singular fibers}\label{skod}
In this section we describe the possible fibers for a minimal elliptic
surface $\pi: X \ra \Ps^1$. The dual graph associated to a (singular, reducible)
curve has a vertex for each irreducible component of the curve and  has an edge between
two vertices if and only if the two corresponding components intersect.

\begin{Thm}[Kodaira] Let $\pi:X \ra \Ps^1$ be an elliptic
  surface. Then the following types of fibers are possible:
\begin{itemize}
\item $I_0$: A smooth elliptic curve.
\item $I_1$: A nodal rational curve. (The dual graph is $\tilde{A}_0$.)
\item $I_\nu, \nu \geq 2$: A $\nu$-gon of smooth rational curves. (The dual graph is $\tilde{A}_{\nu-1}$.)
\item $II$: A cuspidal rational curve. (The dual graph is $\tilde{A}_0$.)
\item $III$: Two rational curves, intersecting in exactly one point
  with multiplicity 2. (The dual graph is $\tilde{A}_1$.)
\item $IV$: Three concurrent lines. (The dual graph is $\tilde{A}_2$.)
\item $I^*_\nu, \nu \geq 0$: The dual graph is of type
  $\tilde{D}_{4+\nu}$.
\item $IV^*, III^*, II^*$: The dual graph is of type $\tilde{E_6},\tilde{E_7},\tilde{E_8}$.
\end{itemize}
\end{Thm}

For more on this see for example \cite[Lecture I]{MiES}, \cite[Appendix C, Theorem
15.2]{Silv} or  \cite[Theorem IV.8.2]{Silv2}.

In Figure~\ref{kodexa}  we give the  dual graph for some of the fiber
types. For a resolution $X\ra Y$ of a fixed double  covering $\varphi':Y\ra\Ps^2$ ramified along $R$ we define we define ``special'' curve as the strict transform  of a component of $\varphi'^{-1}(R)$, and a curve is called 
``ordinary'' otherwise. This notion  depends heavily on our
situation: it gives information on the behavior of the double cover
involution on the fiber components. 

Using our classification of elliptic fibrations one can actually show that $X$ can be obtained in an unique way as a double cover of $\Ps^2$ ramified along six lines up to automorphism, i.e., given two morphisms $\pi_i:X \ra \Ps^2$ of degree 2, ramified along six lines $\ell_j^{(i)}$, with $j=1\dots 6$, there exist an automorphism of $\Ps^2$, mapping $\{\ell_j^{(1)}\}$ to $\{\ell_j^{(2)}\}$.

Moreover, for all fiber types, except $III$ and $I_2$,
one knows which components are special rational curves. In the dual
graphs given here a  vertex is drawn as a full circle  if  the
component is a ordinary rational curve; a vertex is drawn as an open cirlce then
the corresponding component is a special rational curve.

\setlength{\unitlength}{1mm}
\begin{figure}[hbtp]
\begin{picture}(40,30)
\put(6,5){\line(1,0){18}}
\put(5,6){\line(0,1){18}}
\put(25,6){\line(0,1){18}}
\put(6,25){\line(1,0){18}}
\put(5,5){\circle*{3}}
\put(25,25){\circle*{3}}
\put(5,25){\circle{3}}
\put(25,5){\circle{3}}
\put(13,13){{$I_4$}}
\end{picture}
\begin{picture}(70,30)
\put(4,9){\line(1,1){5}}
\put(4,21){\line(1,-1){5}}
\put(12,15){\line(1,0){7}}
\put(21,15){\line(1,0){8}}
\put(31,16){\line(1,1){5}}
\put(31,14){\line(1,-1){5}}
\put(10,15){\circle{3}}
\put(20,15){\circle*{3}}
\put(30,15){\circle{3}}
\put(3,8){\circle*{3}}
\put(3,22){\circle*{3}}
\put(37,8){\circle*{3}}
\put(37,22){\circle*{3}}
\put(18,19){{$I_2^*$}}
\end{picture}
\caption{Examples of dual graphs of singular fibers}\label{kodexa}
\end{figure}
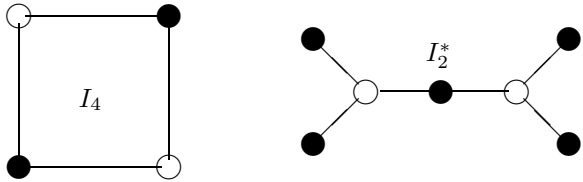

We now give some invariants associated to the singular fibers of an elliptic surface $\pi:X\ra\Ps^1$.
\begin{Def} Let $\pi: X \rightarrow \Ps^1$ be an elliptic surface. Let $P$ be a point of $\Ps^1$. Define $v_P(\Delta_P)$ as the valuation at $P$ of the minimal discriminant of the Weierstrass model, which equals the topological Euler characteristic of $\pi^{-1}(P)$. \end{Def}

\begin{Prop}\label{Noether} Let $\pi: X \rightarrow \Ps^1$ be an elliptic surface, not birational to a product $E\times \Ps^1$, with $E$ an elliptic curve. Then
\[ \sum_{P\in C} v_P(\Delta_P) = 12(p_g(X)+1). \]
\end{Prop}

\begin{proof} This follows from Noether's formula (see \cite[p. 20]{BPV}). The precise reasoning can be found in \cite[Section III.4]{MiES}.
\end{proof}

\begin{Rem}\label{singtabel} If $P$ is a point on $\Ps^1$, such that $\pi^{-1}(P)$ is
  singular then $j(\pi)(P)$ and $v_p(\Delta_p)$ behave as in Table~\ref{kodtab}.
For proofs of these facts see \cite[p. 150]{BPV}, \cite[Theorem IV.8.2]{Silv2} or \cite[Lecture 1]{MiES}.
\end{Rem}

\begin{table}
\[ \begin{array}{|c|c|c|c|}
\hline
\mbox{Kodaira type of fiber over }P & j(\pi)(P) & v_p(\Delta_p) & \mbox{ number of components}\\
\hline
I_0^*                  & \neq \infty & 6 &1\\
I_\nu \;(\nu>0) & \infty &\nu &\nu+1\\
I_\nu^* \;(\nu>0) & \infty&6+\nu  &\nu+5 \\
II     & 0  & 2& 1\\
IV &0& 4& 3\\
IV^* & 0 & 8& 7 \\
II^* & 0 & 10 &9\\
III &1728 & 3& 2\\
III^*              & 1728 & 9&8\\ \hline \end{array}\]
\caption{Invariants of singular fiber types.}\label{kodtab}
\end{table}
\section{Divisors on double covers of $\Ps^2$ ramified along six lines}
\label{sNS} In this section we study the N\'eron-Severi group of
a double cover of $\Ps^2$ ramified along six lines. Most of the
results are probably well known to the experts and are likely to be
found somewhere in the literature.

\begin{Not} Fix six distinct lines $L_i\subset \Ps^2$, such that no
  three of them are concurrent. Denote by $P_{i,j}$  the point of intersection of $L_i$ and $L_j$.

Let $\varphi': Y \ra \Ps^2$ be the double cover ramified along the six
lines $L_i$. Then $Y$ has 15 double points of type $A_1$. Resolving
these points gives a surface $X$, with fifteen exceptional divisors
and a rational map $\varphi: X\ra \Ps^2$. Denote by $\sigma$ the
involution on $X$ induced by the double-cover involution on $Y$ associated to $\varphi'$.

Let $\tilde{\Ps}$ be the blow-up of $\Ps^2$ at the points
$P_{i,j}$. Then  $X/\langle \sigma \rangle =\tilde{\Ps}$, and $\varphi
: X \ra \Ps^2$ factors through $\psi: X \ra \tilde{\Ps}$. Then $\psi$
is a degree 2 cover with branch locus  $\tilde{B}$, the strict
transform of the six lines $L_i$.

 Denote by $\ell_{i,j}\subset X$ the divisor obtained by blowing up
 the point of $Y$ above $P_{i,j} \; (i<j)$. Let $\ell_i$ be such that $2\ell_i$ is the strict transform of $\varphi'^*(L_i)$. 

Let $M^{i,j}_{k,m}$ be the line connecting $P_{i,j}$ and
$P_{k,m}$ ($i<j, k<m, i<k$ and $i,j,k,m$ pairwise distinct). Let
$\mu^{i,j}_{k,m}\subset X$ be the strict transform of $\varphi'^* M^{i,j}_{k,m}$.
\end{Not}

With these notations we have the following intersection results.
\begin{Lem} We have  $\ell_i\cdot \ell_j=-2\delta_{i,j}$, $\ell_{i,j}\ell_{k,m}=-2\delta_{i,k}\delta_{j,m}$ and
  $\ell_i \cdot \ell_{k,m}=\delta_{i,k}+\delta_{i,m}$.\end{Lem}

\begin{Lem} \label{irred} The curve $\mu^{i,j}_{k,m}$ is irreducible if and only if
  $M^{i,j}_{k,m}$ does not intersect any of the 13 points
  $P_{i',j'}$ with $(i',j')\neq (i,j),(k,m)$.\end{Lem}

\begin{proof} 
   The curve $\mu^{i,j}_{k,m}$ is reducible if and only if the strict transform $\tilde{M}$ of
   $M^{i,j}_{k,m}$ on $\tilde{\Ps}$ intersects the branch locus
   $\tilde{B}$ of $\psi : X \ra \tilde{P}$ at every point of intersection
   with even multiplicity. 

  Let $n,n'\in \{1,\dots,6\}\setminus \{i,j,k,m\}$
   such that $n'\neq n$. Let $P$ be the intersection point of
   $M^{i,j}_{k,m}$ and $L_n$. Since the intersection multiplicity of
   $M^{i,j}_{k,m}$ and $B$ at $P$ is even this implies that $P$ is
   also on $L_{n'}$, hence $P=P_{n,n'}$. Conversely, if
   $M^{i,j}_{k,m}$ passes through $P_{n,n'}$ then $\tilde{M}$ and
   $\tilde{B}$ are disjoint, hence the degree morphism
   from  $\mu^{i,j}_{k,m}$ to the rational curve $\varphi'^*(M^{i,j}_{k,m})$ is unramified,
   hence $\mu^{i,j}_{k,m}$ is reducible.
\end{proof}

\begin{Lem}\label{Picnumber} Let $N$ be the subgroup of $NS(X)$
   generated by the $\ell_i$'s, the $\ell_{i,j}$'s and the
   $\mu^{i,j}_{k,m}$'s. Then $N\cong \Z^{16}$ as abelian groups. Moreover $\ell_1$ and the $\ell_{i,j}$'s form a basis for $N\otimes \Q$. \end{Lem}

\begin{proof} First we show that $\rank N \leq 16$. Since $2\ell_i+ \sum_j \ell_{i,j}$ is linearly equivalent
  to the pull-back of a general line $L\subset \Ps^2$, we obtain that
  $\varphi^*L \in N$. Since $\varphi^* L$ is linear equivalent to $\mu^{i,j}_{k,m} +
  \sum a_{k,m}\ell_{k,m}$ for some choice of $a_{k,m}$, we obtain that
the  $\ell_2,\dots \ell_6, \mu^{i,j}_{k,m}$ are contained in \
\[( \Z [\ell_1]\oplus \langle [ \ell_{i,j}] \rangle) \otimes \Q, \]
hence   $\rank(N)\leq 16$.

Since all $\ell_{i,j}$ are disjoint and $\ell_{i,j}^{2}=-2$, they form
a subgroup $N'$ of rank 15 of the N\'eron-Severi group. An easy
computation shows that $\ell_1$ is not linearly equivalent to a divisor
contained in $N' \otimes \Q$, proving that $\rank N=16$. Since the intersection pairing is non-degenerate we have that $N$ is torsion-free.
\end{proof}
\begin{Lem}\label{Action} The action of $\sigma$ on $N$ is trivial. \end{Lem}

\begin{proof} This follows from the fact that for all $i,j$ the automorphism $\sigma|_{\ell_i}$ is the identity 
   and $\sigma$ fixes the curves $\ell_{i,j}$.\end{proof}

\begin{Lem}\label{split} Let $L\subset \Ps^2$ be a rational curve, $L\neq L_i$. Suppose
  $\varphi'^{-1}(L)$ has two components. Then $\rho(X)\geq
  17$.\end{Lem}

\begin{proof} Let $r_i, i=1,2$ be the strict transforms on $X$ of the
  components of $\varphi'^{-1}(L)$. It suffices to show that
  $r_1\not \in N\otimes \Q$. Since $\sigma(r_1)=r_2$ and
  $\sigma$ acts trivial on $N$, we obtain that $r_1$ is linearly
  equivalent to $r_2$. From Lemma~\ref{selfint} it follows that $r_1^2=-2$. We have that any irreducible effective divisor $r\neq r_1$ linear equivalent to $r_1$ satisfies $\# r\cap r_1=r.r_1=-2$ , hence  the only effective
  irreducible divisor linearly equivalent to $r_1$ is $r_1$
  itself. From this it follows that $r_1=r_2$, contradicting our assumption.\end{proof}

\begin{Lem} Suppose that there exists a permutation $\tau \in S_6$ such that $P_{\tau(1),\tau(2)}$, $P_{\tau(3),\tau(4)}$, $P_{\tau(5),\tau(6)}$ are collinear. Then $NS(X)$ has rank at least 17.\end{Lem}

\begin{proof} Using Lemma~\ref{irred} the assumption implies  that at least one of the $\mu^{i,j}_{k,m}$
is reducible. Hence by Lemma~\ref{split} we have  that 
  $\rank NS(X)\geq 17$.\end{proof}

\begin{Ass}\label{assjac} For the rest of this article assume that the six lines in $\Ps^2$ are chosen in such a
  way  that $\rank NS(X)=16$. (Hence the $\mu^{i,j}_{k,m}$ are reduced and irreducible).\end{Ass}

\begin{Rem}\label{rmkdimcal} Choosing 6 distinct lines in $\Ps^2$ gives 4
  moduli. Hence the family $G$ of
$K3$ surfaces that can be obtained as a double cover of six lines in
$\Ps^2$ is 4-dimensional. From the Torelli theorem for $K3$
surfaces (\cite{PSS}) it follows that a general element $X\in G$  satisfies
$\rho(X)\leq 16$.  By Lemma~\ref{Picnumber} any $X\in G$  has
  Picard number at least 16. Hence a  general $X\in G$ satisfies $\rho(X)=16$.
Therefore the above assumption makes sense.\end{Rem}


\begin{Def} Let $C$ be an irreducible curve on $X$ different from the $\ell_i$'s, with $C^2=-2$. Then $C$ is a called an {\em ordinary rational curve}. If $C$ is one of the $\ell_i$'s then $C$ is called a {\em special rational curve}. \end{Def}

A smooth rational curve $C$ satisfies $p_a(C)=0$. From
Lemma~\ref{selfint} it follows that $C^2=-2$, which explains one of the conditions in the above definition.

Let $B=\ell_1+\ell_2+\cdots +\ell_6$.  Since $\psi : X \ra
\tilde{\Ps}$ is ramified along the strict transforms of the $L_i$, we
obtain that the fixed locus  $X^\sigma$
of $\sigma$ equals $B$.

\begin{Lem}\label{intersect} Let $D$ be an ordinary rational
  curve. Then $D\cdot B=2$. \end{Lem}
\begin{proof} Assumption~\ref{assjac} and Lemma~\ref{Picnumber}  imply that the involution $\sigma$ maps
  $D$ to a linear equivalent divisor. Since $D^2=-2$ the only effective
  irreducible divisor linear equivalent to $D$, it $D$ itself. Hence $\sigma$ acts
  non-trivially  on $D\cong \Ps^1$, from which it follows that  there are two fixed points. \end{proof}

\begin{Lem}\label{intersect2} Let $D_1$ and $D_2$ be two ordinary
  rational curves. Then $D_1\cdot D_2\equiv 0 \bmod 2$. \end{Lem}

\begin{proof} The proof of \cite[Lemma 1.6]{Ogui} carries over in this case. \end{proof}

\section{A special result}\label{sNL}

\begin{Prop}\label{NLbnd} Suppose $X$ is a double cover of $\Ps^2$
  ramified along six lines. If the position of the six lines is
  sufficiently general  then for every elliptic fibration $\pi: X \ra
  \Ps^1$  the total number of fibers of $\pi$ of type $II, III$ or $IV$ is at most $\rank(MW(\pi))$. \end{Prop}

\begin{proof}As explained in Remark~\ref{rmkdimcal}  the general
  member of the family of all double covers of $\Ps^2$ ramified along 6
  lines has Picard number 16, hence 
\[ 20=h^{1,1}(X)=4+\rho_{tr}(\pi)+\rank MW(\pi).\]

Using \cite[Proposition 4.6]{NL}
we obtain that
\begin{eqnarray*} 4 &\leq& \dim \{[\psi: X\ra\Ps^1] \in \M_2 | C(\psi)=C(\pi) \}\\ &=& h^{1,1}-\rho_{tr}(\pi)-\#\{\mbox{fibers of type }II, III, IV\} \\&\leq &4+\rank MW(\pi) - \# \{\mbox{fibers of type } II, III, IV \},\end{eqnarray*}
which gives the desired inequality.\end{proof}


\begin{Rem} This proposition can be used to determine the number of
  fibers of type $I_1$ and $I_2$ in several cases of Oguiso's
  classification \cite{Ogui}. However, Oguiso classified all Jacobian elliptic fibrations on the Kummer surface of the product of two non-isogenous elliptic curves, while if one wants to apply Proposition~\ref{NLbnd} one obtains only the classification of Jacobian elliptic fibrations on a Kummer surface of a product of two general elliptic curves.
%
\end{Rem}

\section{Possible singular fibers}\label{sps}
In this section we classify all elliptic fibrations on the double cover of $\Ps^2$ ramified along six fixed lines $\ell_i$ in general position, where general position means that $\rho(X)=16$. 

\begin{Def} By a {\em simple component} $D$ of a fiber $F$ we mean an irreducible component
$D$ of $F$ such that  $D$ occurs with multiplicity one in $F$.\end{Def}

\begin{Prop} \label{posfib}Let $X$ be as before. Let $\pi: X \ra \Ps^1$ be an elliptic
  fibration with a section. Then the Kodaira type of the singular
  fiber is contained in the following list. For each Kodaira type we
  list  the number of components which are special rational curves,
  the number of simple components which are special rational curves, and
  the number of simple components which are ordinary rational curves,
  are contained in the following list:
\[\begin{array}{l|ccc}
\mbox{Type} & \#\mbox{Cmp. special}  & \#\mbox{Simple
cmp. special} & \#\mbox{Simple cmp. ordinary} \\
& \mbox{rational curves} &\mbox{rational curves} & \mbox{rational curves}\\
\hline

I_1   & 0 & 0 & 0 \\
I_2   & 0\mbox{ or }1 & 0\mbox{ or }1& 2\mbox{ or }1 \\ 
I_4   & 2  & 2 & 2\\
I_6  & 3 &3 &3\\
I_8  & 4 &4 &4\\
I_{10} &5 & 5 &5\\
I_0^* & 1 & 0 & 4\\
I_2^* & 2 & 0 & 4\\
I_4^* & 3 & 0 & 4\\
I_6^* & 4 & 0 & 4\\
II & 0 & 0 & 0 \\
III &  0 \mbox{ or } 1 & 0\mbox{ or }1& 2\mbox{ or }1 \\ 
IV^* & 4 & 3 & 0\\
III^* & 3 & 0 & 2\\
II^* & 4 & 0 &1. \\
\end{array}\]\end{Prop}

\begin{proof}[Sketch of proof] 
First of all we prove that no fiber $F$ of type $I_{2k+1}, k>0$
exists. Such a fiber $F$ is a $2k+1$-gon of rational curves. Since two special rational curves
do not intersect, and two ordinary rational curves have even
intersection number (Lemma~\ref{intersect2}), it follows that every ordinary rational curve
intersects two special rational curves, and each special rational
curve intersects two ordinary rational curves. This forces the number
of components in an $n$-gon to be even. Hence $I_{2k+1}$
does not occur. For the same reasons no fiber of type $IV$ or type
$I_{2k+1}^*$ occurs.

Since there at most 6 special nodal curves, no fiber of type $I_{2k},
k>6$ or of type $I^*_{2k}$, $k>5$ occurs.

We prove now that no fiber $F$ of type $I_{12}$ exists. If it would
exists, then this fiber would contain all special rational
curves. Hence the zero section is an ordinary rational curve $Z$. From
Lemma~\ref{intersect} it follows that then $1=Z\cdot F \geq Z \cdot B
= 2$, a contradiction. The non-existence of $I_{10}^*$ follows
similarly. A fiber of type $I^*_8$ has four ordinary rational curves $R_i$
with the property that $R_i$ intersects only one other fiber
component. Moreover, the $R_i$ are the only simple components. Hence
the zero-section $Z$ intersects one of the $R_i$, say $R_1$, and $Z$
has to be a special rational curve, by Lemma~\ref{intersect2}. The
curve $R_i, i\neq 1$ intersect a special rational curve not contained in
the fiber and different from $Z$, hence there are at most 4 special
rational curves contained in $F$. Using Lemma~\ref{intersect2} one
obtains easily that $F$ contains at least 5 special nodal curves, a contradiction.

Let $D$ be a rational curve intersecting three other disjoint rational
curves $D_i, i=1,\dots, 3$.  If $D$ were ordinary, then by
Lemma~\ref{intersect2} the curves $D_i$ would be special. This would imply that
$D\cdot B\geq 3$, contradicting Lemma~\ref{intersect}. Hence $D$ is a
special rational curve. This observation determines in many cases the
number of special components in a singular fiber.
\end{proof}

\section{Possible configurations I}\label{sca}
In this section we study all Jacobian elliptic fibrations $\pi: X \ra
\Ps^1$ having the property that
all special rational curves are fiber components. This section yields the proof of Theorem~\ref{thminfinite}.
\begin{Prop} Let $X$ be as before, in particular $\rank NS(X)=16$. Let $\pi: X \ra \Ps^1$ be a Jacobian elliptic
  fibration. Suppose all special rational curves are contained in the
  fibers of $\pi$. Then one of the following occurs:
\[\begin{array}{lcc}
\mbox{Singular fibers} & & \mbox{Mordell-Weil rank}\\
I_{10} \; I_{2} \; aII\;bI_1& 2a+b=12& 4\\
I_{8} \;I_{4} \; aII\;bI_1& 2a+b=12&  4\\
2I_{6} \; aII\;bI_1& 2a+b=12& 4\\
IV^*\; I_4 \; aII\;bI_1& 2a+b=12&  5\\
\end{array}\]
Conversely, for each case there exist $a,b$ such that these fibration do occur.
\end{Prop}

\begin{proof} Since the zero section is a rational curve and all special
  rational curves are contained in some fibers, we have that   the
  zero section $Z$ is an ordinary rational curve. From
  Lemma~\ref{intersect2} it follows that if  $Z$ intersects a
  reducible fiber, then it intersects in a simple component, which is
  an special nodal curve. From Lemma~\ref{intersect} and the fact that
  special nodal curves are smooth, it follows that
  there are precisely two reducible fibers.
Using Proposition~\ref{posfib} we obtain that the possible reducible fibers are
  $III, IV^*$ and $I_{2k}$, with $1\leq k \leq 5$. 
 
 Since there are six special rational curves, the above possibilities
  are the only ones, except that  the case $I_{10}\; III \; aII\;bI_1$
  with $2a+b=11$ might occur. We prove
  that this one cannot exist:
A fiber of type $I_{10}$ contains precisely five special rational curves, say
  $\ell_1,\dots \ell_5$. Then $\ell_6$ is a component of the
  $III$-fiber. The other component of this fiber is an ordinary rational
  curve $D$ tangent to $\ell_6$, not intersecting  $\ell_i, i=1,\dots,
  5$. Consider $\psi: X\ra \tilde{\Ps}$. Then $\psi(D)$ is a reduced
  curve tangent to $\psi(\ell_6)$ and not intersecting
  $\psi(\ell_i),i\leq 5$, hence intersects the branch locus with even
  multiplicity in each intersection point. Since $\psi(D)$ is a
  rational curve, this implies that $\psi^{-1}(\psi(D))$
  has two components. Hence $\varphi^{-1}(\varphi(D))$ has two
  components. This contradicts Lemma~\ref{split} (note that we assumed
  that $\rho(X)=16$).

The quantity $2a+b$ can be determined using Noether's formula
(Theorem~\ref{Noether}), the Mordell-Weil rank can be
obtained using the Shioda-Tate formula (Theorem~\ref{ST}).

It remains to prove the existence of the remaining four cases.

Let $k\in \{3,4,5 \}$. To prove the existence of  $I_{2k} I_{12-2k}
aII bI_1$,
take
$D=\ell_{1}+\ell_{1,2}+\ell_{2}+\ell_{2,3}+\cdots+\ell_{k}+\ell_{1,k}$.
If $k=3,4$ then $D_1=\ell_{k+1}+\ell_{k+1,k+2}+\cdots
+\ell_{6}+\ell_{k+1,6}$ is an effective divisor with $D_1^2=0$ and
$D.D_1=0$. From Lemma~\ref{fiblem} it follows that the fibration
associated to $|D|$ has $D$ and $D_1$ as fibers. They are of type
$I_{2k}$ and $I_{12-2k}$.

 If $k=5$, then from Lemma~\ref{fiblem} it follows that  $|D|$ defines
 a fibration with a $I_{10}$ fiber, which proves the existence of the
 first case.

To prove the existence of the fibration with the $IV^*$ and $I_4$
fiber, take $D=\ell_1+2 \ell_{1,2} +3\ell_2+2\ell_{2,3}
+\ell_3+2\ell_{2,4}+\ell_4$. Then the fibration associated to $|D|$
has a fiber of type $IV^*$. Yielding the final case.

It remains to prove that the above fibrations are Jacobian. In all
cases one easily shows that $D.\ell_{1,6}=1$, hence $\ell_{1,6}$ is a section.
\end{proof}
%
%

\section{Possible configurations II}\label{scb}In this section we consider fibrations on $X$ such that at least one of the special rational curves is not a component of a singular fiber. 
This section yields the proof of Theorem~\ref{thmfinite}.

Let $\pi:X \ra \Ps^1$ be a Jacobian elliptic fibration, such that at least one of the special rational curves is not contained in a fiber.

\begin{Lem} The group $MW(\pi)$ is finite. \end{Lem}

\begin{proof} The proof of \cite[Lemma 2.4]{Ogui} carries over.\end{proof}

\begin{Lem}\label{sec} Suppose one of the singular fibers of $\pi$ is of type $I_{2k}^*$, $III^*$ or $II^*$. Then all sections are special rational curves. \end{Lem}

\begin{proof} A section intersects a reducible fiber in a simple
  component with multiplicity one. By Proposition~\ref{posfib} we know that all simple
  components in the above mentioned singular fibers are ordinary
  rational curves. From Lemma~\ref{intersect2} we know that two
  ordinary rational curves intersect with even multiplicity, hence
  every section is a special rational curve.
\end{proof}

\begin{Lem} Let $\pi: X \ra \Ps^1$ be as above. Then the only fibers of
  type $I_\nu$ are of type $I_1$ and $I_2$, and no fiber of type
  $I_2$ or type $III$  contains a special rational curve as a
  component. 
\end{Lem}

\begin{proof} From Proposition~\ref{posfib} we know that every fiber
  of type $I_\nu$, $\nu>1$ is of type $I_{2k}$, $k\geq 1$

Without loss of generality we may assume that $\ell_1$
  is not contained in a fiber of $\pi$, hence intersects every
  fiber. Let $F$ be a  singular fiber of type $I_{2k}, k \geq 1$ or of type $III$, containing a
  special rational curve. 

Then $\ell_1$ intersects a reducible fiber in an ordinary
  rational curve, say $D$. If $F$ is of type $I_{2k}$, $k>1$ then $D$ intersects two other
  components, and by Lemma~\ref{intersect2} these components are
  special rational curves. Let $B=\sum \ell_i$. Then $D\cdot B\geq3$, contradicting
  Lemma~\ref{intersect}. Hence it is not possible to have a fiber of
  type $I_{2k}$, $k>1$ containing a special rational curve. By
  Proposition~\ref{posfib} 
  every singular fiber of type $I_{2k}$, $k>1$ contains a special
  nodal curve, hence such a fiber does not occur.

If $k=1$ or the fiber is of type $III$, then  $D$ intersects the other
  component twice, and, by assumption, this component is
  special. Hence $D \cdot B\geq 3$, contradicting Lemma~\ref{intersect}. 
\end{proof}

From now on we study the possibilities for the fibration $\pi$. We distinguish eight cases. In each case we suppose that $\pi$ has a fiber $F$ of a certain given Kodaira
type. In each case  we study which other singular fibers can
occur. Then we prove
the existence of the configuration by giving a divisor $D$ such that
the linear system $|D|$ gives the desired fibration. First, we determine
the fiber types of fibers containing a special rational curve. Observe
that by Proposition~\ref{posfib} all singular fibers, not containing special rational curves, are of type $III, I_2, II$ or $I_1$. Using
the Shioda-Tate formula (Theorem~\ref{ST}) and
Noether's formula (Theorem~\ref{Noether}) we can
determine the quantities $iii+i_2$ and $3iii+2 i_2+2ii+i_1$.

For existence proofs we use 
Proposition~\ref{posfib}. 

\begin{Def} An {\em end-component} $C$ of a fiber $F$ is  a
component $C$ intersecting the support of $\overline{F\setminus C}$ transversally in one point.\end{Def}
In the sequel we use that if an end-component of
a fiber  is a ordinary rational curve, then this component
has to intersect a special rational curve, not contained in the
fiber. This follows immediately from Lemma~\ref{intersect}.

In order to determine the Mordell-Weil group, we use 
Lemma~\ref{sec}. 


\subsection{$II^*$}
In this case four special rational curves are contained in $F$. Using
Proposition~\ref{posfib} it
follows that there are three  end-components $E_1,E_2,E_3$, of which
$E_1$ and $E_2$ are ordinary
rational curves,  and $E_1$ is a simple component. This
means that the zero-section is a special nodal curve, say $D_1$, and
$E_2$ intersects a special nodal curve, not contained in $F$ and
different from $D_1$, say $D_2$. This yields that  $F\cdot D_1=1$,
$F\cdot D_2=3$. All
other fibers are of type $II, III, I_1, I_2$, yielding case
2.1. For example take
$D=\ell_{1,5}+2\ell_1+3\ell_{1,2}+4\ell_2+5\ell_{2,3}+6\ell_3+4\ell_{3,4}+2\ell_4+3\ell_{3,6}$.
Then $\ell_5$ is a section and $\ell_6$ is a trisection.
\subsection{$III^*$}
In this case three special rational curves are contained in $F$. At least two special rational curves have positive intersection number with $F$.

Suppose that two special rational curves are sections, then the third special rational curve is a multisection and all other fibers are of type $II,III, I_1, I_2$.
For example, take $D=\ell_{3,4}+2\ell_3+3\ell_{1,3} +4 \ell_1+2
\ell_{1,5}+3\ell_{1,2}+2\ell_2+\ell_{2,6}$. Then $\ell_4$ and $\ell_6$
are sections. This gives the case 2.2.

If precisely one special rational curve is a section, then there is a special rational which is a multisection, and one special rational curve which is contained in some singular fiber. Since the multisection and the section intersect the fiber  four times, this fiber cannot be of type $III$ or $I_2$, so it is of type $I_0^*$. All other fibers are of type $II, III, I_1, I_2$.
For example take $D=\ell_{3,4}+2\ell_3+3\ell_{1,3} +4 \ell_1+2
\ell_{1,5}+3\ell_{1,2}+2\ell_2+\ell_{2,5}$. This gives case
2.3. In this case $\ell_4$ is a section.

\subsection{$IV^*$}
In this case $F$ contains four special rational curves. The other two
special rational curves do not intersect this fiber, so they are
components of other singular fibers, hence all special rational curves
are components, contradicting our assumptions.
\subsection{$I_6^*$}
In this case $F$ has four special rational curves as  components. One
special rational curve is a section and one a multisection. All other
fibers are of type $III, I_2, II$ and $I_1$. For example take
$D=\ell_{1,5}+\ell_{1,6}+2\ell_1+2\ell_{1,2}+2\ell_2+2\ell_{2,3}+2\ell_3+2\ell_{3,4}+2\ell_4+\ell_{4,5}+\mu^{1,3}_{2,6}$.
The curve $\ell_6$ is a section. This gives case 2.4.
\subsection{$I_4^*$}
In this case $F$ has three special rational curves as  components. Either three or two of the other special curves intersect any fiber. 

If all three special curves intersect any fiber then all other
singular fibers are of type $III$, $I_2$, $II$, $I_1$. For example take
$D=\ell_{1,5}+\ell_{1,4}+2 \ell_1+2 \ell_{1,2}+2 \ell_2+2 \ell_{2,3}+2
\ell_{3}+ \ell_{3,5}+ \ell_{3,6}$. In this case $\ell_4$ and $\ell_6$
are sections. This gives case 2.5.

If two special curves intersect any fiber, then one of the special
curves is again a component of a singular fiber. Since the two other
special rational curves intersect any fiber four times, this fiber cannot
be of type $I_2$ or $III$, so it is a fiber of type $I_0^*$. For
example take $D=\ell_{1,5}+\ell_{1,4}+2 \ell_1+2 \ell_{1,2}+2 \ell_2+2
\ell_{2,3}+2 \ell_{3}+ \ell_{3,5}+ \mu^{1,6}_{2,4}$. In this case,
$\ell_4$ is a section. This gives case 2.6.

\subsection{$I_2^*$}
In this case $F$ fiber has 2 special curves as  components.

If 4 special curves intersect any fiber then all other singular fibers are of type $III$, $I_2$, $II$, $I_1$.
For example take \[D=\ell_{1,3} +
\ell_{1,4}+2\ell_1+2\ell_{1,2}+2\ell_2 + \ell_{2,5} +\ell_{2,6}.\] Then
$\ell_i, i=3,\dots, 6$ are sections. From \cite[Corollary VII.3.3]{MiES} it follows
that $MW(\pi)\not \cong \Z/4\Z$. This gives case 2.7.

If 3 special curves intersect any fiber then there is one special
rational curve $D_1$ not intersecting the fiber and not contained in the
$F$. From this it follows that $D_1$ is a component of a fiber of type
$I_0^*$. All other fibers are of type $II, III, I_1, I_2$. For example
take $D=\ell_{1,3} + \ell_{1,4}+2\ell_1+2\ell_{1,2}+2\ell_2 +
\ell_{2,4} +\ell_{2,5}$. Then $\ell_3$ and $\ell_5$ are sections, This gives case 2.8.

If 2 special curves intersect any fiber then there are two remaining
special rational curves $D_1,D_2$, with $F\cdot D_1=F\cdot D_2=0$. If $D_1$ and
$D_2$ are components of the same fiber, then this fiber is of type
$I_2^*$. For example take $D=\ell_{1,3} +
\ell_{1,4}+2\ell_1+2\ell_{1,2}+2\ell_2 + \ell_{2,4}
+\mu^{1,5}_{3,6}$. Then $\ell_5, \ell_6$ and $\ell_{5,6}$ are
components of another singular fiber, which has to be of type
$I_2^*$. Then $\ell_3$ is a section. This gives case 2.9.

If $D_1$ and $D_2$ are in different fibers then they are both
components of fibers of type $I_0^*$. For example take $D=\ell_{1,3} +
\ell_{1,4}+2\ell_1+2\ell_{1,2}+2\ell_2 + \ell_{2,4} +C'$, with $C'$
the strict transform of $\varphi'^{-1}(C)$, with $C$ the
the conic through $P_{1,3} P_{1,5}P_{2,3} P_{4,6} P_{5,6}$. The curves
$\ell_{3,5},\ell_{3,6},\ell_5,\ell_6$ are components of some singular
fibers. Since $F\cdot \ell_3=1$, they are components of two distinct
fibers. This gives case 2.10.

\subsection{$I_0^*$}
In this case $F$ contains only one special rational curve. Each ordinary
component intersects only one  special rational curve not contained in
$F$. Hence there are at most 4 special curves intersecting $F$.

If there are four special rational curves $D_i$, with $F\cdot D_i>0$, then the
 remaining special rational curve is a component of a fiber, which has to
 be of type $I_0^*$.
 There are still six components of fibers left. The only way to
 arrange them is with 6 $I_2$ fibers. For example take
 $D=2\ell_1+\sum_{k=2}^5\ell_{1,k}$, then
 $D'=2\ell_6+\sum_{k=2}^5\ell_{k,6}$ is another singular fiber. The
 rational curves
 $\ell_k, k=2,\dots, 5$ are sections. From \cite[Corollary VII.3.3]{MiES} it follows
 that $MW(\pi)\not \cong \Z/4\Z$. This gives case 2.11.

If there are three special rational curves $D_i$, with $F\cdot D_i>0$, then the two other special rational curves are components of some singular fiber. If both components are in the same fiber then that fiber is of type $I_2^*$ (which is handled above), otherwise the fibers containing the special rational curves are of type $I_0^*$. So we have in total 3 fibers of type $I_0^*$. All other fibers are of type $II, III, I_1, I_2$. For example take $D=\mu^{2,3}_{5,6}+\ell_{1,4}+\ell_{1,5}+\ell_{1,6}+2\ell_{1}$.
Then none of  $\ell_2, \ell_{2,4}, \ell_{2,5}, \ell_{2,6}, \ell_3,
\ell_{3,4}, \ell_{3,5}, \ell_{3,6}$ intersect a fiber. Hence they are
components of two $I_0^*$ fibers. The curves $\ell_5$ and $\ell_6$ are
sections. This gives case 2.12.

If there are two special rational curves $D_i$, with $F\cdot D_i>0$, then the three other special curves $D'_i$ are components of some singular fibers.

If all $D_i'$'s are contained in the same fiber, then that fiber is either of type $III^*$ (which we already handled) or of type $I_4^*$ (which we also handled above).

If all $D_i'$'s are contained in two singular fibers then we obtain one fiber of type $I_2^*$ and one of $I_0^*$. This case we handled above.

If all $D_i'$'s are contained in three singular fibers then all three fibers are of type $I_0^*$, which is impossible, since then the Picard number of $X$ would be at least 18.

\subsection{Only $I_2$ and $III$}
From the Shioda-Tate formula (Theorem~\ref{ST}) it
follows that there are $\rho(X)-2-\rank MW(\pi)=14$ singular fibers of
type $I_2$ or $III$. From Noether's formula~\ref{Noether}
it follows that then
$24=12p_g(X)+12=\sum v_p(\Delta_p) \geq 2 \cdot 14=28$. A
contradiction. Hence this does not occur.


\begin{thebibliography}{10}

\bibitem{BPV}
W.~Barth, C.~Peters, and A.~Van~de Ven.
\newblock {\em Compact complex surfaces}.
\newblock Springer, 1984.

\bibitem{Har}
R.~Hartshorne.
\newblock {\em Algebraic geometry}, volume~52 of {\em Graduate Texts in
  Mathematics}.
\newblock Springer-Verlag, New York-Heidelberg, 1977.

\bibitem{proefschrift}
R.~Kloosterman.
\newblock {\em Arithmetic and moduli of elliptic surfaces}.
\newblock PhD thesis, University of Groningen, 2005.

\bibitem{NL}
R.~Kloosterman.
\newblock Higher noether-lefschetz loci for elliptic surfaces.
\newblock Preprint available at arxiv:math/AG/0501454, 2005.

\bibitem{MiES}
R.~Miranda.
\newblock {\em The basic theory of elliptic surfaces}.
\newblock Dottorato di Ricerca in Matematica. ETS Editrice, Pisa, 1989.

\bibitem{Ogui}
K.~Oguiso.
\newblock On jacobian fibrations on the {K}ummer surfaces of the product of
  non-isogenous curves.
\newblock {\em Journal of the Mathematical Society of Japan}, 41:651--680,
  1989.

\bibitem{PSS}
I.~I. Pjatecki\u{\i}-{\v{S}}apiro and I.~R. {\v{S}}afarevi{\v{c}}.
\newblock Torelli's theorem for algebraic surfaces of type {${K3}$}.
\newblock {\em Izv. Akad. Nauk SSSR Ser. Mat.}, 35:530--572, 1971.

\bibitem{Sd}
T.~Shioda.
\newblock On the {M}ordell-{W}eil lattices.
\newblock {\em Comment. Math. Univ. St. Paul.}, 39:211--240, 1990.

\bibitem{Silv}
J.H. Silverman.
\newblock {\em The {A}rithmetic of {E}lliptic {C}urves}, volume 106 of {\em
  GTM}.
\newblock Springer-Verlag, New York, 1986.

\bibitem{Silv2}
J.H. Silverman.
\newblock {\em Advanced topics in the arithmetic of elliptic curves}, volume
  151 of {\em GTM}.
\newblock Springer-Verlag, New York, 1994.

\end{thebibliography}

\end{document}